\documentclass[10pt]{amsart}
\usepackage{amsmath,amssymb,amsfonts,amsthm,amsxtra}
\usepackage{cases}
\usepackage{pdfsync}
\usepackage{color}
\usepackage{colordvi}
\usepackage{multicol,longtable}
\usepackage{hyperref}
\usepackage{mathtools}
\usepackage[margin=1.1in]{geometry}
\usepackage{xcolor}
\hypersetup{
    colorlinks,
    linkcolor={red!50!black},
    citecolor={blue!50!black},
    urlcolor={blue!80!black}
}

\makeatletter
\newcommand\appendix@section[1]{%
  \refstepcounter{section}%
  \orig@section*{Appendix \@Alph\c@section: #1}%
  \addcontentsline{toc}{section}{Appendix \@Alph\c@section: #1}%
}

\newtheorem{theorem}{Theorem}
\newtheorem*{theorem*}{Theorem}
\newtheorem{lemma}{Lemma}
\newtheorem*{lemma*}{Lemma}
\newtheorem{corollary}{Corollary}
\newtheorem*{corollary*}{Corollary}
\newtheorem{proposition}{Proposition}
\newtheorem*{proposition*}{Proposition}
\newtheorem{remark}{Remark}
\newtheorem*{remark*}{Remark}

\newcommand{\thm}[1]{\begin{theorem}#1\end{theorem}}

\newcommand{\lem}[1]{\begin{lemma}#1\end{lemma}}

\newcommand{\prop}[1]{\begin{proposition}#1\end{proposition}}

\newcommand{\eq}[1]{\begin{equation}#1\end{equation}}
\newcommand{\eqs}[1]{\begin{equation*}#1\end{equation*}}
\newcommand{\alig}[1]{\begin{align}#1\end{align}}
\newcommand{\aligs}[1]{\begin{align*}#1\end{align*}}
\newcommand{\nn}{\nonumber}
\newcommand{\pf}[1]{\begin{proof}#1\end{proof}}

\newcommand{\mbb}[1]{\mathbb{#1}}
\newcommand{\mcal}[1]{\mathcal{#1}}

\newcommand{\pt}{\hspace*{0.5pt}}

\newcommand{\mmod}{\mathrm{mod}}
\newcommand{\ts}{\textstyle}
\newcommand{\ds}{\displaystyle}
\newcommand{\sss}{\scriptscriptstyle}
\newcommand{\ip}[2]{\langle #1,#2\rangle}
\newcommand{\dd}{\hspace*{0.1em}\mathrm{d}}
\newcommand{\pd}{\partial}
\newcommand{\f}[2]{\frac{#1}{#2}}
\newcommand{\el}{\left}
\newcommand{\er}{\right}
\newcommand{\re}{\mathrm{Re}\pt}
\newcommand{\im}{\mathrm{Im}}
\newcommand{\res}[1]{\underset{#1}{\,\,\mathrm{res}\,\,}}
\newcommand{\rarrow}{\rightarrow}

\newcommand{\sym}{\mathrm{sym}}

\newcommand{\ep}{\varepsilon}

\newcommand{\hf}{\f{1}{2}}
\newcommand{\thf}{{\textstyle\f{1}{2}}}
\newcommand{\qt}{\f{1}{4}}
\newcommand{\tqt}{{\textstyle\f{1}{4}}}

\newcommand{\tf}[2]{{\textstyle\frac{#1}{#2}}}

\newcommand{\heq}{\hspace*{3.1ex}}

\newcommand{\qand}{\quad\mbox{and}\quad}

\newcommand{\subf}{_{\hspace*{-1pt}f}}
\newcommand{\subF}{_{\hspace*{-1pt}F}}

\newcommand{\intt}{\int_{-\infty}^\infty}

\newcommand{\nb}{\newblock}
\begin{document}

\title{Nonvanishing of central $L$-values of Maass forms}

\author[S. Liu]{Shenhui Liu}

\address{Department of Mathematics, The Ohio State University\\ 231 W 18th Avenue\\
Columbus, OH 43210}
\email{liu.2076@osu.edu}

\begin{abstract}
With the method of moments and the mollification method, we study the central $L$-values of GL(2) Maass forms of weight $0$ and level $1$ and establish a positive-proportional nonvanishing result of such values in the aspect of large spectral parameter in short intervals, which is qualitatively optimal in view of Weyl's law. As an application of this result and a formula of Katok--Sarnak, we give a nonvanishing result on the first Fourier coefficients of Maass forms of weight $\hf$ and level $4$ in the Kohnen plus space.
\end{abstract}



\maketitle
\tableofcontents
\section{Introduction}\label{intro}
Nonvanshing of central $L$-values and their derivatives of automorphic forms is an important research topic, due to the connection between such values and various aspects of mathematics, such as arithmetic geometry, spectral deformation theory, and analytic number theory. The combination of the method of moments and the mollification method, initiated by Iwaniec--Sarnak \cite{IwaniecSarnak2000}, has been a very fruitful approach in yielding positive-proportional nonvanishing results on central $L$-values and their derivatives in a family of automorphic forms (see, e.g., \cite{IwaniecSarnak1999}, \cite{VanderKam1999}, \cite{KowalskiMichel1999}, \cite{KowalskiMichel2000}, \cite{KowalskiMichelVanderKam2000}, \cite{Soundararajan2000}, \cite{KowalskiMichelVanderKam2002}, \cite{Blomer2008}, \cite{Khan2010}, \cite{Luo2015}, \cite{Liu2016}, and others). Along this direction we address the case of GL(2) Maass forms. Specifically, we study the (mollified) moments of the $L$-functions of the Hecke--Maass forms of weight $0$ and level $1$ at the central point of the critical strip, and establish a positive-proportional nonvanishing result of such values in short intervals of the spectral parameters (Theorem \ref{mainthm}). As an application, this result and a formula of Katok--Sarnak (see (\ref{KS})) imply a strong nonvanishing result (Theorem \ref{weight1/2}) of the first Fourier coefficient of Maass forms in the Kohnen plus space of weight $\hf$ and level $4$.  \\


Let $S_0(1)$ be the space of Maass cusp forms of weight $0$ and level $1$ and pick an orthonormal basis $\{u_j\}$ of Hecke--Maass forms of $S_0(1)$, where each $u_j$ has Laplace eigenvalue $\qt+t_j^2$ ($t_j\geq0$). (See \S\,\ref{sec: review} for a brief review of Maass forms.) In the rest of this work we always let $T>0$ be a large parameter and assume $T^\eta<M<T(\log T)^{-1}$ with a fixed small $0<\eta<1$. Our main result is the following  
\thm{\label{mainthm}We have
$$
\#\el\{t_j\mid|t_j-T|\leq M, L(\thf,u_j)>0\er\}\gg TM.
$$
}
By Weyl's law (see \cite{Venkov1982} and \cite{Hejhal1983})
$$
N(T):=\#\{j\mid t_j\leq T\}=\f{1}{12}T^2-\f{1}{2\pi}T\log T+c_0T+O(T(\log T)^{-1}),
$$
we have
$$
N(T+M)-N(T-M)=\f{1}{3}TM+O\!\el(T\er),
$$
i.e., there are $\asymp TM$ many $t_j$'s in the interval $[T-M,T+M]$. Hence Theorem \ref{mainthm} implies that for Hecke-Maass forms in the basis $\{u_j\}$ with spectral parameter $t_j\in[T-M,T+M]$, there are positive proportion of them with nonvanishing central $L$-values. This is analogous to Luo's nonvanishing result \cite{Luo2015} for central $L$-values of holomorphic cusp forms for $\Gamma_0(1)$ of large weight, which is our main motivation. It is worth mentioning that Xu \cite{Xu2014} obtains a positive-proportion nonvanishing result of the $L(\hf+it_j,u_j)$ for $t_j$ in short intervals, using mollifiers and moments but with different treatment.

In view of the author's work \cite{Liu2016} on central $L$-derivative values of holomorphic cusp forms for $\Gamma_0(1)$ of large weight, one expects a similar nonvanishing result for $L'(\hf,u_j)$ for odd Hecke-Maass eigenforms $u_j$ ($\ep_j=-1$). A possible approach to prove this, say, is to adapt Motohashi's formula (Lemma \ref{moto}) to treat a twisted moment of $L'(\hf,u_j)^2$ and apply the mollification analysis in \cite{Liu2016}. \\

Now let $S_0(4)$ be the space of Maass cusp forms of weight $\hf$ and level $4$ and denote by $S_\hf^+(4)$ the Kohnen plus space. Pick an orthonormal basis $\{F_j\}$ of Hecke--Maass forms in $S_\hf^+(4)$, where each $F_j$ has Laplace eigenvalue $\qt+t^2_{\hspace*{-1pt}F_j}$ with $t_{\hspace*{-1pt}F_j}\geq0$. (Again see \S\,\ref{sec: review} for a review.) As an application of Theorem \ref{mainthm} and a formula of Katok--Sanark $(\ref{KS})$, and the fact that every weight $\hf$ Maass form lifts to a weight 0 Maass form, we have the following
\thm{\label{weight1/2}
For the Hecke--Maass forms in the basis $\{F_j\}$ with $t_{\hspace*{-1pt}F_j}\in[T-M,T+M]$, there are positive proportion of them whose first Fourier coefficient $b_{\hspace*{-1pt}F_j}(1)\neq0$.
}
\vspace*{1ex}
%
%
%
%

In the following we outline the structure of the paper and give the proof of Theorem \ref{mainthm} and some comments. We approach the nonvanishing problem in Theorem \ref{mainthm} via the study of the harmonic moments
\aligs{
\sum_j \f{L(\hf,u_j)^kM_j^k}{L(1,\sym^2u_j)}h_0(t_j)\quad (k=1,2)\qand \sum_j \f{h_0(t_j)}{L(1,\sym^2u_j)^2}.
}
Here the test function $h_0(t)$ is given by
\alig{\label{hh_0}
h_0(t)=T^{-2}h(t),
}
where 
$$
h(t)=(t^2+\tqt)\pt\omega(t)\qand\omega(t):=\omega_{T,M}(t)=\ds e^{-\el(\!\f{t-T}{M}\!\er)^2}+e^{-\el(\!\f{t+T}{M}\!\er)^2};
$$
and $M_j$ ($j\geq1$) are mollifiers defined in (\ref{mollifier_defn}). We remark that $h_0(t)$ gives a more natural counting than $h(t)$ but in the actual computation we use $h(t)$ in place of $h_0(t)$ to avoid writing the factor $T^{-2}$ everywhere. One reason for including the extra factor $t^2+\qt$ in $h(t)$ is that Motohashi's formula (Lemma \ref{moto}), which we use to treat the second moment, requires that $h(\pm \f{i}{2})=0$. \vspace*{1ex}


For completeness we record the following asymptotic formulas for the unmollified moments with power-saving, which seem not to have been stated in the literature.
\prop{\label{unmollifiedmoments}
We have
\aligs{
 &\sum_j \f{L(\hf,u_j)}{L(1,\sym^2u_j)}h_0(t_j)+\f{1}{2\pi}\intt\f{|\zeta(\hf+it)|^2}{|\zeta(1+2it)|^2}h_0(t)\dd t=\f{2}{\pi^{3/2}}TM+O(T^{-\hf+\ep}M)\\
 \intertext{and}
 &\sum_j \f{L(\hf,u_j)^2}{L(1,\sym^2u_j)}h_0(t_j)+\f{1}{2\pi}\intt\f{|\zeta(\hf+it)|^4}{|\zeta(1+2it)|^2}h_0(t)\dd t\\
 &\hspace*{10.5em} =\f{2}{\pi^{3/2}}\big(TM\log T+(\gamma-\log2\pi)TM\big)+O(T^{-2}M^3\log T),
}
where $\gamma$ is the Euler constant.
}

The power-saving in the above indicates that there is room to insert mollifiers \`{a} la Selberg to kill the extra $\log T$ in the second moment, i.e., to bring the mollified moments to comparable size as in Lemma \ref{mollifiedmoments}, whose proof constitutes the major part of this investigation. The proof of Proposition \ref{unmollifiedmoments} can be viewed as a simplified version of that for Lemma \ref{mollifiedmoments}. For example, the second asymptotic formula in Proposition \ref{unmollifiedmoments} follows from Motohashi's formula (Lemma \ref{moto}) for $n=1$ and the estimates (\ref{psi_p1}) and (\ref{psi_m1}).
We remark that closely related to the asymptotics in Proposition \ref{unmollifiedmoments} are the following upper bounds for the unmollified and unweighted moments
\aligs{
\sum_{T-M\leq t_j\leq T+M}\f{L(\hf,u_j)}{L(1,\sym^2u_j)}&\ll TM\\
\intertext{and}
\sum_{T-M\leq t_j\leq T+M}\f{L(\hf,u_j)^2}{L(1,\sym^2u_j)}&\ll TM\log T,
}
due to Ivi\'{c}--Jutila \cite{IvicJutila2003} and Motohashi \cite{Motohashi1992}, respectively. \\

Next we explain the use of the mollified moments and prove Theorem \ref{mainthm}. After some preparation in \S\,\ref{sec: preparation}, we establish the following estimates for the mollified moments in \S\S\,\ref{sec: mollified_1st_moment}--\ref{sec: mollified_2nd_moment}.
\lem{\label{mollifiedmoments}Let $\delta$ be the number which appears in the definition of $M_j$ (see $(\ref{mollifier_defn})$). If $0<\delta<\f{3}{10}$ we have
\aligs{
\sum_j \f{L(\hf,u_j)M_j}{L(1,\sym^2u_j)}h_0(t_j)&=\f{1}{\pi^{3/2}}TM+o(TM)\\
\intertext{and if $0<\delta<\f{1}{4}$ we have}
\sum_j \f{L(\hf,u_j)^2M_j^2}{L(1,\sym^2u_j)}h_0(t_j)&\ll TM.
}
}
For the mollified first moment we apply an approximate functional equation (Lemma \ref{ApproxFE}) for $L(\hf,u_j)$ and the Kuznetsov trace formula over even forms (Lemma \ref{KTFeven}). (Note that $L(\hf,u_j)=0$ for odd forms $u_j$.) The treatment of the off-diagonal sum $\mcal{O}^+$ involving the $J$-Bessel function $J_{2it}(x)$ is inspired by Li's work \cite{Li2009,Li2011}. While for the off-diagonal sum $\mcal{O}^-$ involving the $K$-Bessel function $K_{2it}(x)$, we split the $c$-sum of Kloosterman sums $S(m,n;c)$ into two ranges, treat small $c$ by Li's idea, and for large $c$ do a stationary phase analysis using an asymptotic formula of $K_{2it}(x)$.

For the mollified second moment, we employ Motohashi's formula (Lemma \ref{moto}) at the outset, instead of using an approximate functional equation for $L(\hf,u_j)^2$. The benefit is that the right-side of Motohashi's formula does not involve any Kloosterman sums or Bessel functions, but only shifted sums of the divisor function and certain functions $\Psi^\pm$ for which Motohashi's work \cite{Motohashi1992,Motohashi1997} and Ivi\'{c}'s work \cite{Ivic2001} provide convenient resources. On the other hand, Luo's work \cite{Luo2015} also reduces the expected high load of analysis for the mollified second moment, since Luo's successful mollification analysis can be applied directly right after we apply Motohashi's formula.

Now we give the deeper reason for using Motohashi's formula. If one would proceed with an approximate functional equation for $L(\hf,u_j)^2$ and Kuznetsov over even forms, one then wishes to perform analysis analogous to holomorphic modular form cases as in \cite{LauTsang2005,Luo2015,Liu2016}, namely, to extract information from the off-diagonal terms resulting from Kuznetsov over even forms by using properties of Estermann zeta-functions. But this would not be easy since the Mellin--Barnes representation of $J_{2it}(x)$ gives very narrow room for contour shifting. And in fact, this is not necessary, for in the derivation of Motohashi's formula (\cite[\S\,3.3]{Motohashi1997}) one already uses analysis involving Estermann zeta-functions, and more importantly the outset of the derivation gives the advantage of getting rid of the ``cumbersome'' $J$-Bessel term, which is inevitable if one uses Kuznetsov over even forms (see also the penultimate paragraph on p.\,113 of \cite{Motohashi1997}).  \\


In Section \ref{sec: neg_1st_moment_symmsq}, we prove the following upper bound, which is a short-interval version of \cite[Lemma 5]{Luo2001}.
\lem{\label{negative_1st_moment_symmsq}
\aligs{
\sum_{j}\f{h_0(t_j)}{L(1,\sym^2u_j)^2}\ll TM.
}
}

Finally we are ready to prove Theorem \ref{mainthm}. By Lemma \ref{mollifiedmoments}, Lemma \ref{negative_1st_moment_symmsq} and H\"{o}lder's inequality, we have
\aligs{
TM&\ll\sum_{j}\f{L(\hf,u_j)M_j}{L(1,\sym^2u_j)}h_0(t_j)\\
&\ll\bigg(\sum_{L(\hf,u_j)\neq0}h_0(t_j)\bigg)^{\!\qt}\bigg(\sum_{j}\f{h_0(t_j)}{L(1,\sym^2u_j)^2}\bigg)^{\!\qt}
\bigg(\sum_j \f{L(\thf,u_j)^2M_j^2}{L(1,\sym^2u_j)}h_0(t_j)\bigg)^{\!\hf}\\
&\ll\bigg(\sum_{L(\hf,u_j)\neq0}h_0(t_j)\bigg)^{\!\qt}(TM)^{\f{3}{4}},
}
Hence we have
\aligs{
TM&\ll \sum_{L(\hf,u_j)\neq0}h_0(t_j)\\
&\ll \sum_{\substack{|t_j-T|\leq M\log T\\L(\hf,u_j)\neq0}}\f{t_j^2}{T^2}e^{-\el(\!\f{t_j-T}{M}\!\er)^2}\ll \sum_{\substack{|t_j-T|\leq M\log T\\ L(\hf,u_j)\neq0}}e^{-\el(\!\f{t_j-T}{M}\!\er)^2}
}
and thus
\alig{\label{mainthm_aux1}
TM\ll\sum_{L(\hf,u_j)\neq0}e^{-\el(\!\f{t_j-T}{M}\!\er)^2}.
}
Next we follow Luo \cite{Luo2001} to remove the weight. By partial summation we see that for any fixed $A>0$
\aligs{
\sum_{AM\leq t_j-T\leq 2AM}e^{-\el(\!\f{t_j-T}{M}\!\er)^2}\ll TM\int_A^{2A}e^{-t^2}\dd t.
}
Then applying this inequality to $2^kAM\leq t_j-T\leq 2^{k+1}AM$ and summing over $k$, we get
\alig{\label{mainthm_aux2}
\sum_{t_j\geq T+AM}e^{-\el(\!\f{t_j-T}{M}\!\er)^2}\ll TM \int_A^\infty e^{-t^2}\dd t\ll TM e^{-A^2}.
}
A similar argument shows that
\alig{\label{mainthm_aux3}
\sum_{t_j\leq T-AM}e^{-\el(\!\f{t_j-T}{M}\!\er)^2}\ll TM e^{-A^2}
}
With a sufficiently large $A$, (\ref{mainthm_aux1})$\sim$(\ref{mainthm_aux3}) imply that
\aligs{
TM\ll\sum_{\substack{|t_j-T|\leq AM\\ L(\hf,u_j)\neq0}}e^{-\el(\!\f{t_j-T}{M}\!\er)^2}.
}
Replacing $M$ by $M/A$ in the above yields
\aligs{
TM\ll\sum_{\substack{|t_j-T|\leq M \\ L(\hf,u_j)\neq0}}e^{-\el(\!\f{t_j-T}{M}\!\er)^2}\ll\sum_{\substack{|t_j-T|\leq M\\L(\hf,u_j)\neq0}}1
}
and Theorem \ref{mainthm} follows. \\

\noindent{\bf Acknowledgements.} The author thanks Professor Wenzhi Luo for suggesting this project and for his constant support. The author is grateful to Professor Gergely Harcos for his careful reading and for pointing out several typos.

%
%
%
%
%
%
%
\section{Preparation}\label{sec: preparation}
\subsection{A review of Maass forms of weight $0$ and weight $\hf$}{\label{sec: review}}
Consider the group $\Gamma_0(1)=SL(2,\mbb{Z})$ and its Hecke congruence subgroup $\Gamma_0(4)$, which act on the upper half-plane $\mbb{H}=\{x+iy\mid x\in\mbb{R}\mbox{ and }y>0\}$ by linear fractional transformation, with fundamental domains $\mcal{D}_1$ and $\mcal{D}_4$, respectively. Let $L^2(\mcal{D}_j)$ be the space of functions on $\mbb{H}$ which are square-integrable on $\mcal{D}_j$ with respect to the invariant measure $y^{-2}\dd x\dd y$. Define the Laplace operators
$$
\Delta_{k}=-y^2\!\el(\f{\pd^2}{\pd x^2}+\f{\pd^2}{\pd y^2}\er)\!+iky\f{\pd}{\pd x},\quad k=0\mbox{ or }\hf.
$$
The space $S_0(1)$ of Maass cusp forms of weight 0 and level 1 is the set
\aligs{
S_0(1)=\el\{f\in L^2(\mcal{D}_1) \el|\begin{aligned}
           &f(\gamma z)=f,\ \forall\gamma\in\Gamma_0(1),\ f\mbox{ is cuspidal}\\
           &\Delta_0f=(\tqt+t\subf^2)f\mbox{ for some }t\subf\geq0
          \end{aligned}\er.\er\}.
}
The cuspidality here and below means that the zeroth Fourier coefficient of a form vanishes at all cusps of the relevant fundamental domain. Each $f\in S_0(1)$ has a Fourier expansion at $i\infty$
$$
f(z)=\sum_{n\neq0}a\subf(n)W_{0,it\subf}(4\pi|n|y)e(nx).
$$
Here $e(z)$ denotes $e^{2\pi iz}$; $W_{\mu,\nu}$ denotes the Whittaker function (see \cite[Chapter 7]{MOS1966}), which has a specialization
$$
W_{0,\nu}(4\pi y)=(4y)^\hf K_\nu(2\pi y)
$$
where $K_\nu(z)$ is the usual $K$-Bessel function. The Fourier coefficients $a\subf(n)$ satisfy
$$
a\subf(n)=\ep\subf\pt a\subf(-n)\mbox{ for }n\in\mbb{Z}
$$
where $\ep\subf=1$ or $-1$, according to which we call a form $f$ \textit{even} or \textit{odd}.
For $n\geq1$ the Hecke operator $T_n$ is defined by
$$
T_nf(z)=\f{1}{\sqrt{n}}\sum_{ad=n}\sum_{b\,\mmod\,d}f\!\el(\f{az+b}{d}\er).
$$
If $f\in S_0(1)$ is an eigenfunction of all $T_n$ with eigenvalues $\lambda_f(n)$, we call $f$ a \textit{Hecke--Maass form} and note that $a_f(n)\sqrt{n}=a_f(1)\lambda_f(n)$ ($n\geq1$).
For later use, we fix an orthonormal basis $\{u_j\}$ of $S_0(1)$ consisting of Hecke--Maass forms $u_j$ of Laplace eigenvalues $\qt+t_j^2$ ($t_j\geq0$) and Hecke eigenvalues $\lambda_j(n)$.

 To any Hecke--Maass form $f$ we associate its $L$-function
$$
L(s,f)=\sum_{n\geq1}\f{\lambda_f(n)}{n^s}\quad\mbox{for }\re(s)>1
$$
which admits analytic continuation to the whole complex plane and satisfies the functional equation
$$
\Lambda(s,f):=L_\infty(s,t\subf)L(s,f)=\ep\subf\Lambda(1-s,f)
$$
where
$$
L_\infty(s,t)=\pi^{-s}\Gamma\!\el(\!\f{s+it}{2}\er)\Gamma\!\el(\!\f{s-it}{2}\er).
$$
One expects $L(\hf,f)\neq0$ for many forms $f$ with $\ep\subf=1$, while $L(\hf,f)$ is necessarily $0$ when $\ep\subf=-1$ due to the functional equation.
Related is the symmetric square $L$-function
\aligs{
L(s,\sym^2f)=\sum_{n\geq1}\f{\lambda_f(n^2)}{n^s}\quad\mbox{ for }\re(s)> 1,
}
which also has entire continuation to the whole complex plane. \\

The space $S_\hf(4)$ of Maass cusp forms of weight $\hf$ and level 4 is the set
\aligs{
S_\hf(4)=\el\{F\in L^2(\mcal{D}_4) \el|\begin{aligned}
           &F(\gamma z)=J(\gamma,z)f(z),\ \forall\gamma\in\Gamma_0(1),\ F\mbox{ is cuspidal}\\
           &\Delta_\hf F=(\tqt+t\subF^2)F\mbox{ for some }F\geq0
          \end{aligned}\er.\er\},
}
where the automorphy factor $J(\gamma,z)=\theta(\gamma z)/\theta(z)$ with $\theta(z)=y^\qt\sum_{n\in\mbb{Z}} e(n^2z)$.
Each $F\in S_\hf(4)$ has a Fourier expansion at $i\infty$
$$
F(z)=\sum_{n\neq0}b\subF(n)W_{\qt\mathrm{sgn}(n),it\subF}(4\pi|n|y)e(nx).
$$
Define Hecke operators $T_p^2$ for all primes $p\neq2$
$$
T_{p^2}F(z)=\sum_{n\neq0}\el\{p\, b\subF(np^2)+\el(\f{n}{p}\er)\!\f{1}{\sqrt{p}}\,b\subF(n)+\f{1}{p}\,b_F\!\el(\f{n}{p^2}\er)\er\}W_{\qt\mathrm{sgn}(n),it\subF}(4\pi|n|y)e(nx),
$$
where $b\subF(r)\neq0$ if $r\not\in\mbb{Z}$ and $(\f{n}{p})$ denotes the Legendre symbol. Define an operator $L:S_\hf(4)\rarrow S_\hf(4)$ by
$$
LF(z)=\qt e\!\el(\f{\pi}{8}\er)\!\el(\f{z}{|z|}\er)^{\!-\hf}\sum_{\nu\,\mmod\,4}F\!\el(\f{4\nu z-1}{16z}\er).
$$
Then $L$ is self-adjoint, commutes with $\Delta_\hf$ and all $T_{p^2}$, and satisfies $(L-1)(L+\hf)=0$ (see \cite[Proposition 1.4]{KatokSarnak1993}). The Kohnen plus space $S^+_\hf(4)\subseteq S_\hf(4)$ is the eigenspace of $L$ with eigenvalue $1$ and
$$
F\in S_\hf(4)\mbox{ lies in }S^+_\hf(4)\quad\mbox{if and only if}\quad b\subF(n)=0\mbox{ for }n\equiv2,3\ (\mmod\ 4).
$$
Then we can find an orthonormal basis $\{F_j\}$ of $S^+_\hf(4)$ consisting of common eigenfunctions of all $T_{p^2}$ ($p\neq2$).
For $F\in S_\hf(4)$ we define its Shimura lift ${\rm Sh}F$ by
$$
{\rm Sh}F(z)=\sum_{n\neq0}a_{{\rm Sh}F}(n)W_{0,2it\subF}(4\pi|n|y)e(nx),
$$
where
$$
a_{{\rm Sh}F}(n)=\sum_{n=\ell m}\f{\sqrt{|m|}}{\ell}b\subF(\ell^2).
$$
Then ${\rm Sh}\circ T_{p^2}=T_p\circ{\rm Sh}$; if $F\in S^+_\hf(4)$ then ${\rm Sh}F\in S_0(1)$ with $\Delta_0$-eigenvalue $\qt+(2t\subF)^2$; if $F\in S^+_\hf(4)$ is a common eigenfunction of $T_{p^2}$ ($p\neq2$), then $T_{p^2}F=\lambda_{{\rm Sh}F}(p)F$; $b\subF(1)\neq0$ if and only if ${\rm Sh}F\neq0$. (See \cite[Proposition 4.1]{KatokSarnak1993}.)
Then we have the following \vspace*{1ex}

\noindent{\bf Katok--Sarnak formula} (\cite[(0.19)]{KatokSarnak1993}){\bf.} For a normalized Hecke--Maass form $f\in S_0(1)$ ($a_f(1)=1$)
\alig{\label{KS}
\f{\Lambda(\hf,f)}{\ip{f}{f}}=12\sqrt{\pi}\sum_{{\rm Sh}F_j=b_{\hspace*{-1pt}F_j}\!(1)f}|b_{\hspace*{-1pt}F_j}(1)|^2.
}
Here in the sum we have ${\rm Sh}F_j=b_{\hspace*{-1pt}F_j}\!(1)f$, which is different from \cite[(19)]{KatokSarnak1993} where the nonzero Shimura lifts are arithmetically normalized. An immediate consequence of this formula is the positivity of $L(\hf,f)$.
We comment that Baruch--Mao \cite{BaruchMao2010} shows that the Kohnen plus space $S_\hf^+(f,4)$ for an individual normalized Hecke--Maass form $f\in S_0(1)$, given by
$$
S_\hf^+(f,4)=\big\{F\in S_\hf^+(4) \mid \Delta_\hf F=\tqt(1+t^2\subf)F,\ T_{p^2}F=\lambda_f(p)F\mbox{ for all }p>2\big\},
$$
is one-dimensional. Then by the Katok--Sarnak formula and Baruch--Mao's result, if $L(\hf,f)>0$, which happens ``frequently'' according to Theorem \ref{mainthm}, then some single $F_j$ generates $S_\hf^+(f,4)$ and the sum on the right-hand side of (\ref{KS}) consists of only one summand. 
\vspace*{1ex}
\subsection{Analytic tools}
In the following we introduce the tools required for the study of the relevant harmonic moments.
\subsubsection{Approximate functional equation}
For even $u_j$ in the eigenbasis of $S_0(1)$, we need an approximate functional equation to represent $L(\hf,u_j)$, whose proof is standard as in \cite[Theorem 5.3]{IwaniecKowalski2004}.
\lem{\label{ApproxFE}
\eqs{L({\textstyle{\hf}},u_j)=2\sum_{m\geq1}\f{\lambda_j(m)}{\sqrt{m}}U(m,t_j), }
where
\eqs{U(y,t)=\f{1}{2\pi i}\int_{(A)}y^{-u}G(u)\gamma(u,t)\dd u,}
with any fixed $A>0$, $G(u)=u^{-1}e^{-u^4}$, and
\eqs{
\gamma(u,t)=\f{L_\infty(\hf+u,t)}{L_\infty(\hf,t)}
=\pi^{-u}\f{\Gamma\Big(\f{\hf+u-it}{2}\Big)}{\Gamma\Big(\f{\hf-it}{2}\Big)}\f{\Gamma\Big(\f{\hf+u+it}{2}\Big)}{\Gamma\Big(\f{\hf+it}{2}\Big)}.
}
}
An easy consequence of the functional equation of the Riemann zeta-function is that
\eq{\label{ApproxFE_Eisenstein}
|\zeta(\thf+it)|^2=2\sum_{m\geq1}\f{d_{it}(m)}{\sqrt{m}}U(m,t).
}
In view of Barnes's formula (see Proposition \ref{BarnesError} in Appendix \ref{sec: BarnesFormula}),
we have for fixed $u$, fixed $\sigma\geq0$, and large $t$
\alig{\label{estGamma}
\gamma(u,t-i\sigma)=\Big(\f{|t|}{2\pi}\Big)^{u}e^{O(P(|u|+\sigma)|t|^{-1})},
}
where $P(x)$ is a cubic polynomial with positive coefficients. Hence for any $A>0$ and $\sigma<A+\hf$, the function $U(y,t)$ is holomorphic in the strip $-\sigma\leq\im(t)\leq0$, and
\eq{\label{estU}
U(y,t-i\sigma)\ll_{A,\sigma}\begin{cases}y^{-A},&\quad |t|\leq1\vspace*{1ex}\\|t|^Ay^{-A},&\quad |t|\geq1.\end{cases}
}
\subsubsection{Kuznetsov trace formulas}
In our notation the Kuznetsov traces formulas are as follows.
\lem{\label{KTF}Let $h(t)$ be an even function which is holomorphic in $|\im(t)|\leq\hf$ with $h(t)\ll(1+|t|)^{-2-\ep}$ with some $\ep>0$. Then for integers $m,n\geq1$
\aligs{
&2\sum_{j}\f{\lambda_j(m)\lambda_j(\pm n)}{L(1,\sym^2u_j)}h(t_j)+\f{1}{\pi}\intt\f{d_{it}(m)d_{it}(n)}{|\zeta(1+2it)|^2}h(t)\dd t\\
&=\delta_{m,n}H+\sum_{c>0}\f{S(m,\pm n;c)}{c}H^\pm\!\!\el(\!\f{4\pi \sqrt{mn}}{c}\er).
}
Here $d_{s}(n)=\sum_{n=ab}(a/b)^{s}=d_{-s}(n)$,
\eqs{H=\f{1}{\pi^2}\intt t\pt \tanh(\pi t)h(t)\dd t, }
\eqs{H^+(x)=\f{2i}{\pi}\intt \f{J_{2it}(x)}{\cosh(\pi t)}t\pt h(t)\dd t, }
\eqs{H^-(x)=\f{4}{\pi^2}\intt t\sinh(\pi t)K_{2it}(x)h(t)\dd t, }
and $J_\nu(z)$ and $K_\nu(z)$ are the usual Bessel functions.
}
This is a restatement of \cite[Theorem 9.3]{Iwaniec2002} or \cite[Theorem 2.2 and 2.4]{Motohashi1997}, in view of the relation
$$
\f{|2a_f(1)|^2}{\cosh(\pi t\subf)\ip{f}{f}}=\f{2}{L(1,\sym^2f)}
$$
for any Hecke-Maass form $f\in S_0(1)$. As a consequence of Lemma \ref{KTF}, we have the Kuznetsov trace formulas over even forms:
\lem{\label{KTFeven}Let $h$ be as in the previous lemma. Then for integers $m,n\geq1$
\aligs{
&2\sum_{\ep_j=1}\f{\lambda_j(m)\lambda_j(n)}{L(1,\sym^2u_j)}h(t_j)+\f{1}{\pi}\intt\f{d_{it}(m)d_{it}(n)}{|\zeta(1+2it)|^2}h(t)\dd t\\
\nn &=\hf\delta_{m,n}H+\hf\sum_{c>0}\f{1}{c}\!\el\{S(m,n;c)H^+\!\!\el(\!\f{4\pi \sqrt{mn}}{c}\er)+S(m,-n;c)H^-\!\!\el(\!\f{4\pi \sqrt{mn}}{c}\er)\!\er\}.
}
}
%
\subsubsection{Motohashi's formula}
To treat the second moment, we employ a formula of Motohashi. For any even entire function $h(t)$ such that $h(\pm\thf i)=0$ and
$$
h(t)\ll e^{-c|t|^2}
$$
for some fixed $c>0$ in any fixed horizontal strip, define
\aligs{
\widehat{h}(s)&=\intt \f{\Gamma(s+it)}{\Gamma(1-s+it)}t\pt h(t)\dd t,\\
\Psi^+(x;h)&=\int_{(\beta)}\Gamma(\thf-s)^2\tan(\pi s)\widehat{h}(s)x^s\dd s,\\
\Psi^-(x;h)&=\int_{(\beta)}\Gamma(\thf-s)^2\f{\widehat{h}(s)}{\cos(\pi s)}x^s\dd s,
}
where $-\f{3}{2}<\beta<\hf$. A restatement of Motohashi's formula {\cite[Lemma 3.8]{Motohashi1997}} in our context is as follows.
\lem{\label{moto}For the test function $h(t)$ as in the last paragraph, we have
$$
\mcal{H}(n;h):=\sum_j \f{L(\thf,u_j)^2\lambda_j(n)}{L(1,\sym^2u_j)}h(t_j)=\sum_{\nu=1}^7\mcal{H}_\nu(n;h),
$$
where
\aligs{
\mcal{H}_1(n;h)&=-\f{i}{\pi^3}\el\{(\gamma-\log(2\pi \sqrt{n})(\widehat{h})'(\thf)+\qt(\widehat{h})''(\thf)\er\}\f{d(n)}{\sqrt{n}},\\
\mcal{H}_2(n;h)&=\f{1}{2\pi^3}\sum_{m\geq1}\f{d(m)d(m+n)}{\sqrt{m}}\Psi^+\Big(\f{m}{n};h\Big),\\
\mcal{H}_3(n;h)&=\f{1}{2\pi^3}\sum_{m\geq1}\f{d(m)d(m+n)}{\sqrt{m+n}}\Psi^-\Big(1+\f{m}{n};h\Big),\\
\mcal{H}_4(n;h)&=\f{1}{2\pi^3}\sum_{m=1}^{n-1}\f{d(m)d(n-m)}{\sqrt{m}}\Psi^-\Big(\f{m}{n};h\Big),\\
\mcal{H}_5(n;h)&=-\f{1}{4\pi^3}\f{d(n)}{\sqrt{n}}\Psi^-(1;h),\\
\mcal{H}_6(n;h)&=-\f{6i}{\pi^2}d_\hf(n)\,h'(-\thf i), \\
\mcal{H}_7(n;h)&=-\f{1}{2\pi}\intt\f{|\zeta(\hf+it)|^4}{|\zeta(1+2it)|^2}d_{it}(t)h(t)\dd t.
}
}

%
\subsubsection{Mollifiers}\label{mollifiers}
For convenience, we define as in Luo's work \cite{Luo2015} the mollifier $M_j$ for $u_j$ by
\eq{\label{mollifier_defn}
M_j=\sum_{n\geq1}\f{a_n\mu(n)\lambda_j(n)}{\sqrt{n}}, }
where
\eq{a_n=\begin{cases}
            \ds{\hf}\f{\log^2(T^{2\delta}/n)-\log^2(T^{\delta}/n)}{\log(T^{\delta})}, &1\leq n\leq T^{\delta}, \\
            \ds{\hf}\f{\log^2(T^{2\delta}/n)}{\log(T^{\delta})}, &T^{\delta}\leq n\leq T^{2\delta},\\
            0, &n\geq T^{2\delta}
          \end{cases}}
for some $\delta>0$. It is easy to see that $0\leq a_n\leq 2\log T$. Also the discontinuous integral
$$
\f{1}{2\pi i}\int_{(2)}\f{y^s}{s^3}\dd s=\begin{cases}\f{1}{2}\log^2(y)&\mbox{ if }y\geq1, \\
0&\mbox{ if }0<y\leq1, \end{cases}
$$
gives the analytic form of $a_n$
\eq{\label{intforan}a_n=\f{1}{2\pi i}\int_{(2)}\f{\big(\f{\xi^2}{n}\big)^s-\big(\f{\xi}{n}\big)^s}{s^3}\f{\dd s}{\log\xi}}
with $\xi=T^\delta$.
%
%
%
%
%
%
%
\section{The mollified first moment}\label{sec: mollified_1st_moment}
In this section we prove the asymptotic formula for the mollified first moment in Lemma \ref{mollifiedmoments}. Note that we use the test function $h(t)$ instead of $h_0(t)$ in the derivation (see \ref{hh_0}). By Lemma \ref{ApproxFE} and the definition of $M_j$, we have
\aligs{
\sum_{\ep_j=1}\f{L({\hf},u_j)M_j}{L(1,\sym^2u_j)}h(t_j)
=\sum_{n\geq1}\f{a_n\mu(n)}{\sqrt{n}}\sum_{m\geq1}\f{1}{\sqrt{m}}\bigg(2\sum_{\ep_j=1}\f{\lambda_j(m)\lambda_j(n)}{L(1,\sym^2u_j)}h_m(t_j)\bigg),
}
where
\eqs{h_m(t)=h(t)U(m,t).}
Since $h_m(t)$ satisfies the conditions in Lemma \ref{KTFeven}, the above becomes
\eq{\label{Kuz1}\sum_{\ep_j=1}\f{L({\hf},u_j)M_j}{L(1,\sym^2u_j)}h(t_j) =\mcal{D}-\mcal{C}+\mcal{O}^++\mcal{O}^-,}
where
\aligs{
\mcal{D}&=\hf\sum_{n\geq1}\f{a_n\mu(n)}{n}H_n,\\
\mcal{C}&=\f{1}{\pi}\sum_{n,m}\f{a_n\mu(n)}{\sqrt{nm}}
  \intt\f{d_{it}(m)d_{it}(n)}{|\zeta(1+2it)|^2}h_m(t)\dd t,\\
\mcal{O}^+&=\hf\sum_{n,m}\f{a_n\mu(n)}{\sqrt{nm}}
  \sum_{c>0}\f{S(m,n;c)}{c}H_m^{+}\!\el(\!\f{4\pi \sqrt{mn}}{c}\er),\\
\mcal{O}^-&=\hf\sum_{n,m}\f{a_n\mu(n)}{\sqrt{nm}}
  \sum_{c>0}\f{S(m,-n;c)}{c}H_m^{-}\!\el(\!\f{4\pi \sqrt{mn}}{c}\er).
}
In the following we analyze the above terms on the right-hand side of (\ref{Kuz1}).
%
\subsection{Diagonal contribution $\mcal{D}$}
We claim that for $0<\delta<\hf$
\alig{\label{D_asymp}
\mcal{D}=\f{1}{\pi^{3/2}}T^3M+o(T^3M).
}
By (\ref{intforan}) and the definition of $U(y,t)$,
\aligs{
\mcal{D}&=\f{1}{\pi^2}\int_0^{\infty}t\tanh(\pi t)h(t)\f{1}{2\pi i}\int_{(2)}\f{\xi^{2w}-\xi^w}{w^3\log\xi}\f{1}{2\pi i}\int_{(3)}\f{G(u)\gamma(u,t)}{\zeta(1+w+u)}\dd u\dd w\dd t.
}
Moving the $w$-integral to $\re(w)=-\hf+\ep$ for a small $\ep>0$, we pick up a simple pole at $u=0$ with residue $\zeta(1+w)^{-1}$ and have
\aligs{
\mcal{D}&=\f{1}{\pi^2}\int_0^{\infty}t\tanh(\pi t)h(t)\bigg\{\f{1}{2\pi i}\int_{(2)}\f{\xi^{2w}-\xi^w}{w^3\log\xi}\f{\dd w}{\zeta(1+w)}\\
&\phantom{aaaaaaa}+\f{1}{(2\pi i)^2}\int_{(-\hf+\ep)}\int_{(2)}\f{\xi^{2w}-\xi^w}{w^3\log\xi}\f{G(u)
\gamma(u,t)}{\zeta(1+w+u)}\dd w\dd u\bigg\}\dd t\\
&=\f{1}{\pi^2}\int_0^{\infty}t\tanh(\pi t)h(t)\el\{1+\f{1}{2\pi i}\int_{C_\ep}\f{\xi^{2w}-\xi^w}{w^3\log\xi}\f{\dd w}{\zeta(1+w)}\er.\\
&\phantom{aaaaaaa}+\f{1}{(2\pi i)^2}\int_{(-\hf+\ep)}\int_{(2)}\f{\xi^{2w}-\xi^w}{w^3\log\xi}\f{G(u)
\gamma(u,t)}{\zeta(1+w+u)}\dd w\dd u\bigg\}\dd t\\
\label{D1aux}&=\f{1}{\pi^2}\int_0^{\infty}t\tanh(\pi t)h(t)\bigg\{1+O_\ep\bigg(\f{1}{\log\xi}\bigg)\\
&\phantom{aaaaaaa}+\f{1}{(2\pi i)^2}\int_{(-\hf+\ep)}\int_{(2)}\f{\xi^{2w}-\xi^w}{w^3\log\xi}\f{G(u)
\gamma(u,t)}{\zeta(1+w+u)}\dd w\dd u\bigg\}\dd t,
}
where $C_\ep$ denotes the contour
\eq{\label{contour}
C_\ep=\{it\mid|t|\geq\ep\}\cup\{\ep e^{i\theta}\mid\tf{\pi}{2}\leq\theta\leq\tf{3\pi}{2}\}
}
which starts from $-i\infty$.

It is easy to compute that
\aligs{
\int_0^{\infty}t\tanh(\pi t)h(t)\dd t=\sqrt{\pi}\, T^3M+O(TM^3).
}
Thus we are left with
\aligs{
&\int_0^{\infty}t\tanh(\pi t)h(t)\int_{(-\hf+\ep)}\int_{(2)}\f{\xi^{2w}-\xi^w}{w^3\log\xi}\f{G(u)\gamma(u,t)}{\zeta(1+w+u)}\dd w
\dd u\dd t \\
&=\int_{(-\hf+\ep)}G(u)\int_{(\hf)}\f{\xi^{2w}-\xi^w}{w^3\zeta(1+w+u)\log\xi}\int_0^{\infty}t\tanh(\pi t)h(t)\gamma(u,t)\dd t\dd w\dd u.
}
By (\ref{estGamma}) and considering $t$ in and outside of $[T-M\log T, T+M\log T]$, we see that
$$
\int_0^{\infty}t\tanh(\pi t)h(t)\gamma(u,t)\dd t\ll T^{\f{5}{2}+\ep}M\pt e^{O(P(|u|))}
$$
So the triple integral in the above is
\aligs{
&\ll\f{\xi\,T^{\f{5}{2}+\ep}M}{\log\xi}
\int_{(-{\hf}+\ep)}|G(u)|e^{O(P(|u|))}\int_{(\hf)}\f{|\dd w|}{|w|^3|\zeta(1+w+u)|}|\dd u|\\
&\ll T^{\f{5}{2}+\delta+\ep}M(\log T)^{-1}.
}
Then the claimed asymptotic formula (\ref{D_asymp}) holds for $\delta<\hf$.
%
\subsection{Continuous spectrum part $\mcal{C}$}
By (\ref{ApproxFE_Eisenstein}), as well as that $a_n\ll\log T$ and $a_n=0$ for $n\geq\xi^2$, we have
\aligs{
\mcal{C}&=\f{1}{\pi}\sum_{n\geq1}\f{a_n\mu(n)}{\sqrt{n}}\int_0^\infty\f{h(t)d_{it}(n)}{|\zeta(1+2it)|^2}\bigg(2\sum_{m\geq1}\f{d_{it}(m)}{\sqrt{m}}U(m,t)\bigg)\dd t\\
&=\f{1}{\pi}\sum_{n\geq1}\f{a_n\mu(n)}{\sqrt{n}}\int_0^\infty\f{|\zeta(\hf+it)|^2}{|\zeta(1+2it)|^2}h(t)d_{it}(n)\dd t\\
&\ll\log T\sum_{n\leq\xi^2}\f{d(n)}{\sqrt{n}}\int_0^\infty\f{|\zeta(\hf+it)|^2}{|\zeta(1+2it)|^2}h(t)\dd t\\
&\ll T^{\delta}(\log T)^4\int_{T-M\log T}^{T+M\log T}t^{2+\f{1}{3}}\omega(t)\dd t\\
&\ll T^{2+\f{1}{3}+\delta}M(\log T)^4,
}
which is $o(T^3M)$ upon letting $\delta<\f{2}{3}$. Here we used the classical bounds
$$\zeta(1+it)^{-1}\ll\log |t|\quad(|t|\geq1)\qand \zeta(\thf+it)\ll |t|^\f{1}{6}.$$
%

\subsection{Off-diagonal sum $\mcal{O}^+$}
In the following we show that $\mcal{O}^+$ is negligible. Here and in the sequel a quantity being negligible means that its size is $O_A(T^{-A})$ for any $A>0$. We start with $H_m^+\Big(\!\f{4\pi \sqrt{mn}}{c}\Big)$ and abuse the notation $X=\f{4\pi \sqrt{mn}}{c}$ for convenience. Note that $J_\nu(z)$ is entire in $\nu$ for fixed $z\neq0$ and by the integral representation (see \cite[3.3(5)]{Watson1944})
\eqs{
J_\nu(z)=\f{2(z/2)^\nu}{\sqrt{\pi}\,\Gamma(\hf+\nu)}\int_0^{\f{\pi}{2}}\sin^{2\nu}(\theta)\cos(z\cos(\theta))\dd\theta, \quad\re(\nu)>-\hf
}
and Stirling's formula we have for $x>0$
\alig{\label{estJBessel}
J_{\sigma+it}(x)\ll_\sigma\begin{cases}x^\sigma, &\quad|t|\leq1,\\x^\sigma |t|^{-\sigma}e^{\pi|t|/2},&\quad|t|\geq1.\end{cases}
}
Note that $\cosh(\pi t)$ has simple zeros at $t=i(\hf-k)$ for integers $k$. Let $K>0$ be an integer and $A>K+\hf$, both to be chosen later. By shifting the $t$-integral to $\im(t)=-K$ for a positive integer $K$, we have
\alig{
H_m^{+}(X)&=\f{2i}{\pi}\intt\f{J_{2it}(X)}{\cosh(\pi t)}t\pt h_m(t)\dd t\\
\nn&=4\sideset{}{'}\sum_{k=2}^K\res{t=i(\hf-k)}\f{J_{2it}(X)}{\cosh(\pi t)}t\pt h_m(t)\\
\nn&\heq+\f{2i}{\pi}\int_{-\infty}^\infty \f{J_{2K+2it}(X)}{\cosh(\pi (t-iK))}(t-iK)h(t-iK)U(m,t-iK)\dd t\\
\nn&=:R_m^+(X)+I_m^+(X).
}
Here the notation $\sum'$ means that at most one of the summand is replaced by zero, since there is at most one $k\geq2$ such that $J_{2k-1}(X)=0$ due to the fact that no two of the functions $J_n(z)$ ($n=0,1,2,\cdots$) have any common strictly positive zeros (see \cite[15.28]{Watson1944}).
First, the residue part becomes
\aligs{
R_m^+(X)&=4\sideset{}{'}\sum_{k=2}^K\lim_{t\rarrow i(\hf-k)}\f{t-i(\hf-k)}{\cosh(\pi t)} J_{2it}(X)\pt t\pt h_m(t)\\
&=\sideset{}{'}\sum_{k=1}^K c_k J_{2k-1}(X)\omega(i(\thf-k))U(m,i(\thf-k)),
}
where
$$
c_k=4\pi^{-1}(-1)^k(k-\thf)(k^2-k).
$$
Here we recall that $\omega(t)=e^{-\el(\!\f{t-T}{M}\!\er)^2}+e^{-\el(\!\f{t+T}{M}\!\er)^2}$. Then with $|\omega(i(\hf-k))|\leq e^{-\f{T^2}{M^2}}e^{\f{K^2}{M^2}}$, (\ref{estJBessel}), and (\ref{estU}), we see that
\eq{\label{R_plus}
R_m^+(X)\ll e^{-\f{T^2}{M^2}}\f{n^{K-\hf}}{m^{A-K+\hf}\,c}.
}
In the following we omit the dependence on $A$ and $K$ of the implied constants.

By
\eqs{
\omega (t-iK)\ll \omega (t)\qand |\cosh(\pi(t-iK))|=\cosh(\pi t)\gg e^{\pi |t|}
}
and (\ref{estJBessel}),
we have
\alig{
\label{I_plus}I_m^+(X)
&\ll \f{n^K}{m^{A-K}c^{2K}}\int_{|t|\leq1}\omega(t)\dd t+\f{n^K}{m^{A-K}c^{2K}}\int_{|t|\geq1}|t|^{A-2K+3}\omega(t)\dd t\\
\nn&\ll T^{A-2K+3}M\f{n^K}{m^{A-K}c^{2K}}.
}
Then (\ref{R_plus}) and (\ref{I_plus}) imply the bound
\eqs{
H_m^+(X)\ll e^{-\f{T^2}{M^2}}\f{n^{K-\hf}}{m^{A-K+\hf}\,c}+T^{A-2K+3}M\,\f{n^K}{m^{A-K}\,c^{2K}},
}
This, together with Weil's bound on Kloosterman sums, yields
\alig{
\nn\mcal{O}^+&=\sum_{n\geq1}\f{a_n\mu(n)}{\sqrt{n}}\sum_{m\geq1}\f{1}{\sqrt{m}}\sum_{c>0}\f{S(m,n;c)}{c}H_m^{+}\!\el(\!\f{4\pi \sqrt{mn}}{c}\er)\\
\nn&\ll e^{-\f{T^2}{M^2}}\log T\sum_{n\leq\xi^2}n^{K-\hf}\sum_{m\geq1}\f{1}{m^{A-K+\hf}}\sum_{c\geq1}\f{d(c)}{c^{3/2}}\\
\nn&\quad\quad\quad+T^{A-2K+3}M\log T\sum_{n\leq\xi^2}n^{K-\hf}\sum_{m\geq1}\f{1}{m^{A-K+1}}\sum_{c\geq1}\f{d(c)}{c^{2K}}\\
\nn&\ll T^{(A-K)+3+\delta+(2\delta-1)K}M\log T,
}
which is negligible upon taking $\delta<\hf$, sufficiently large $K$ and suitable $A$ with $A-K>\hf$.
%
\subsection{Off-diagonal sum $\mcal{O}^-$}We write
\aligs{
\mcal{O}^-=\mcal{O}^-_1+\mcal{O}^-_2,
}
by splitting the $c$-sum into two ranges: $c\geq\sqrt{mn}$ ($X\leq4\pi$) and $c<\sqrt{mn}$ ($X>4\pi$). Here recall the notation $X=\f{4\pi\sqrt{mn}}{c}$.

Case $\mcal{O}^-_1$ ($c\geq\sqrt{mn}$). We start with the identity (\cite[3.7(6)]{Watson1944})
\eqs{K_\nu(z)=\f{\pi}{2}\f{I_{-\nu}(z)-I_{\nu}(z)}{\sin(\pi\nu)}, }
where the $I$-Bessel function $I_\nu(z)$ is entire for fixed $z\neq0$ and has integral representation (\cite[3.71(9)]{Watson1944})
$$
I_\nu(z)=\f{(z/2)^\nu}{\sqrt{\pi}\,\Gamma(\hf+\nu)}\int_0^\pi e^{z\cos(\theta)}(\sin\theta)^{2\nu}\dd \theta,\quad\re(\nu)>-\hf.
$$
Thus for $x>0$
$$
I_{\sigma+it}(x)\ll_\sigma\begin{cases}x^\sigma e^x, &\quad|t|\leq1,\\x^\sigma e^x |t|^{-\sigma}e^{\pi|t|/2},&\quad|t|\geq1.\end{cases}
$$
By $\sin(i2z)=i\sinh(2z)=2i\sinh(z)\cosh(z)$, we have
\aligs{
 H_m^{-}(X)&=\f{4}{\pi^2}\f{\pi}{2}\intt t\sinh(\pi t)h_m(t)\f{I_{-2it}(X)-I_{2it}(X)}{2i\sinh(\pi t)\cosh(\pi t)}\dd t\\
\nn&=\f{i}{\pi}\intt \f{I_{2it}(X)}{\cosh(\pi t)}t\pt h_m(t)\dd t-\f{i}{\pi}\intt \f{I_{-2it}(X)}{\cosh(\pi t)}t\pt h_m(t)\dd t\\
\nn&=\f{2i}{\pi}\intt \f{I_{2it}(X)}{\cosh(\pi t)}t\pt h_m(t)\dd t.
}
Then by a very similar argument as the treatment for $\mcal{O}^{+}$, we see that $\mcal{O}^-_1$ is also negligible in size. \\

Case $\mcal{O}^-_2$ ($c<\sqrt{mn}$). We write
\eqs{H_m^{-}(X)=\f{8}{\pi^2}\int_0^\infty t \sinh(\pi t)K_{2it}(X)h_m(t)\dd t=H_{m,1}^{-}(X)+H_{m,2}^{-}(X)+H_{m,3}^{-}(X), }
by splitting the $t$-integral into $\int_0^{T-M\log T}$, $\int_{T-M\log T}^{T+M\log T}$, and $\int_{T+M\log T}^\infty$. By $K_{2it}(X)\ll e^{-\pi t}t^{-\f{1}{3}}$ (Proposition \ref{Kit_estimate} in Appendix \ref{appendixB}), we have
\eqs{
H_{m,1}^{-}(X)+H_{m,3}^{-}(X)\ll_{A} T^{A+4-\log T}\f{1}{m^A}.
}
It follows that for fixed $A>1$, the contribution of
\aligs{
\label{offdiagtwo1steq2}&\sum_{n,m}\f{a_n\mu(n)}{\sqrt{nm}}\sum_{c<\sqrt{mn}}c^{-1}S(-m,n;c)(H_{m,1}^{-}(X)+H_{m,3}^{-}(X))
}
is negligible.

To achieve
\eq{\label{Knegative_moment}
\sum_{n,m}\f{a_n\mu(n)}{\sqrt{nm}}\sum_{c<\sqrt{mn}}c^{-1}S(-m,n;c)H_{m,2}^-(X)=o(T^3M),
}
we split the sum in $m$ into two ranges: $m\geq T^{1+\ep}$ and $m<T^{1+\ep}$.

In view of the range of $t$ in $H_{m,2}^-(X)$, if $m\geq T^{1+\ep}$, then for sufficiently large $A$
\eqs{U(m,t)\ll_A \f{t^A}{m^A}\ll_{A,\ep}\f{1}{m^2}.}
This and the estimate $K_{2it}(X)\ll e^{-\pi t}t^{-\f{1}{3}}$ imply
$$
H_{m,2}^-(X)\ll T^{3-\f{1}{3}}M\f{1}{m}.
$$
Hence we have
\aligs{
&\sum_{n\leq\xi^2}\f{a_n\mu(n)}{\sqrt{n}}\sum_{m\geq T^{1+\ep}}\f{1}{\sqrt{m}}\sum_{c<\sqrt{mn}}c^{-1}S(-m,n;c)H_{m,2}^-(X)\\
&\ll T^{3-\f{1}{3}}M\log T\sum_{n\leq\xi^2}1\sum_{m\geq T^{1+\ep}}\f{1}{m^{2+1/2}}\sum_{c<\sqrt{mn}}\f{d(c)}{\sqrt{c}}\\
&\ll T^{3-\f{1}{3}}M(\log T)^2\sum_{n\leq\xi^2}n^\qt\sum_{m\geq T^{1+\ep}}\f{1}{m^{2+1/4}}\\
&\ll T^{3+\f{5}{2}\delta-\f{1}{3}-\f{5}{4}}M
}
which is $o(T^3M)$ provided $\delta<\f{19}{30}$.

Now we deal with the case when $m< T^{1+\ep}$. In view of the range of $t$ in $H_{m,2}^-(X)$, we have $X=o(t)$ and the asymptotic formula
(see \cite[p.\,142]{MOS1966})
\eqs{
K_{2it}(X)=\f{\sqrt{2\pi}\,e^{-\pi t}}{(4t^2-X^2)^{\f{1}{4}}}\el\{\sin(\phi_X(t))+O(t^{-1})\er\},
}
where
$$
\phi_X(t)=\f{\pi}{4}+2t\cosh^{-1}\!\Big(\f{2t}{X}\Big)-\sqrt{4t^2-X^2}.
$$
The error term in \cite[p.\,142]{MOS1966} is only $O(X^{-1})$ but the error term $O(t^{-1})$ follows from the power series expansion for $K_{2it}(x)$ through that of $I_{\pm 2it}(x)$. For $T-M\log T\leq t\leq T+M\log T$, we have
$$\sinh(\pi t)e^{-\pi t}=\hf+O(e^{-\pi T}),$$
$$
\f{1}{(4t^2-X^2)^{\f{1}{4}}}=\f{1}{\sqrt{2t}}+O\big(X^2 t^{-5/2}\big),
$$
and
\eqs{
U(m,t)=\f{1}{2\pi i}\int_{(\hf),|\im(u)|\leq T^\ep}\f{t^{u}}{(2\pi m)^u}G(u)\dd u+O_A(T^{-A}).
}
We only need to estimate
\alig{
\label{offdiagtwo1steq4}&\int_{T-M\log T}^{T+M\log T}\bigg(\int_{(\hf),|\im(u)|\leq T^\ep}t^{u+\f{5}{2}}\f{G(u)}{m^u}\dd u\bigg)
e^{-\el(\!\f{t-T}{M}\!\er)^2}\sin(\phi_X(t))\dd t\\
\nn=&\int_{(\hf),|\im(u)|\leq T^\ep}\f{G(u)}{m^u}\bigg(\int_{T-M\log T}^{T+M\log T}t^{u+\f{5}{2}}e^{-\el(\!\f{t-T}{M}\!\er)^2}\sin(\phi_X(t))\dd t\bigg)\dd u.
}
since it is easy to see that the contribution from other parts to $\mcal{O}^-$ through $H_{m,2}^-(X)$ is $o(T^3M)$.
To handle the $t$-integral, we work with
$$
I(X,T)=\int_{T-M\log T}^{T+M\log T}f(t)e^{ig(t)}\dd t,
$$
where
$$
f(t)=t^3e^{-\el(\!\f{t-T}{M}\!\er)^2}\qand g(t)=\im(u)\log t\pm\phi_X(t).
$$
Notice that
\eqs{g'(t)=\im(u)t^{-1}\pm2\log\!\bigg(\!\f{2t+\sqrt{4t^2-X^2}}{X}\bigg)\asymp\pm\log T}
and
$$g''(t)=-\im(u)t^{-2}\pm4(4t^2-X^2)^{-\hf}\asymp T^{-1}. $$
By Fa\`{a} di Bruno's formula for high derivatives of composite functions (see e.g. \cite{Warren2002})
\eqs{\f{\dd^n}{\dd t^n}\el(e^{-\el(\!\f{t-T}{M}\!\er)^2}\er)
=e^{-\el(\!\f{t-T}{M}\!\er)^2}\sum_{m_1+2m_2=n}\f{(-1)^{m_1+m_2}}{m_1!m_2!}\el(\f{2(t-T)}{M}\er)^{\!m_1}\!\f{n!}{M^n}\ll_n M^{-n}(\log T)^n, }
so that
$$
f^{(n)}(t)=\sum_{k=n-3}^n{n \choose k}(t^3)^{(n-k)}\f{\dd^{k}}{\dd t^{k}}\el(e^{-\el(\!\f{t-T}{M}\!\er)^2}\er)\ll_n T^3M^{-n}(\log T)^n
$$

Repeated integration by parts gives
\aligs{
I(X,T)&=-i\int_{T-M\log T}^{T+M\log T}\el(\f{f'(t)}{g'(t)}-\f{f(t)g''(t)}{(g'(t))^2}\er)e^{ig(t)}\dd t+O(T^{3-\log T})\\
&=-i\int_{T-M\log T}^{T+M\log T}\f{f'(t)}{g'(t)}e^{ig(t)}\dd t+O(T^2M(\log T)^{-2})\\
&=\int_{T-M\log T}^{T+M\log T}\el(\f{f''(t)}{(g'(t))^2}-\f{2f'(t)g''(t)}{(g'(t))^3}\er)e^{ig(t)}\dd t+O(T^2M(\log T)^{-2})\\
&=\int_{T-M\log T}^{T+M\log T}\f{f''(t)}{(g'(t))^2}e^{ig(t)}\dd t+O(T^2M(\log T)^{-2})\\
&\cdots\\
&=c_n\int_{T-M\log T}^{T+M\log T}\f{f^{(n)}(t)}{(g'(t))^n}e^{ig(t)}\dd t+O_n(T^2M(\log T)^{-2})\\
}
where $|c_n|=1$. By the estimates for $f^{(n)}(t)$ and $g'(t)$, and by $M>T^\eta$, we take sufficiently large $n$ to obtain
$$
I(X,T)\ll T^2M(\log T)^{-2}.
$$
Hence the integral in (18) is
$$
\ll T^2M(\log T)^{-2}\f{1}{\sqrt{m}}
$$
and its contribution to $\mcal{O}^-$ through $H_{m,2}^-(X)$ is
%
\aligs{
&\ll T^2M(\log T)^{-1}\sum_{n\leq \xi^2}\f{1}{\sqrt{n}}\sum_{m<T^{1+\ep}}\f{1}{m}\sum_{c<\sqrt{mn}}\f{|S(-m,n;c)|}{c}\\
&\ll T^2M(\log T)^{-1}\sum_{n\leq \xi^2}1\sum_{m<T^{1+\ep}}\f{1}{m}\sum_{c<\sqrt{mn}}\f{d(c)}{\sqrt{c}}\\
&\ll T^2M\sum_{n\leq \xi^2}n^\qt \sum_{m<T^{1+\ep}}\f{1}{m^{\f{3}{4}-\ep}}\\
&\ll T^{2+\qt+2\ep+\f{5}{2}\delta}M
}
which is $o(T^3M)$ provided $\delta<\f{3}{10}$. That is, we have shown the bound (\ref{Knegative_moment}). Now the claimed asymptotic formula for the first mollified moment in Lemma \ref{mollifiedmoments} follows.  %
\section{The mollified second moment}\label{sec: mollified_2nd_moment}
In this section we establish the upper bound of the mollified second moment in Lemma \ref{mollifiedmoments}, or the equivalent upper bound (see (\ref{hh_0}))
$$
\sum_j \f{L(\hf,u_j)^2M_j^2}{L(1,\sym^2u_j)}h(t)\ll T^3M.
$$
By the Hecke relation we have
\alig{
M^2_j=\sum_{r\geq1}\f{1}{r}\sum_{n\geq1}\f{\lambda_j(n)}{\sqrt{n}}A_{r,n},
}
where
$$
A_{r,n}:=\sum_{n=n_1n_2}a_{rn_1}\mu(rn_1)a_{rn_2}\mu(rn_2).
$$
Then we apply Lemma \ref{moto} to get
\aligs{
&\sum_j \f{L(\hf,u_j)^2M_j^2}{L(1,\sym^2u_j)}h(t)=\sum_{r\geq1}\f{1}{r}\sum_{n\geq1}\f{A_{r,n}}{\sqrt{n}}\mcal{H}(n;h)\\
&=\sum_{\nu=1}^7\sum_{r\geq1}\f{1}{r}\sum_{n\geq1}\f{A_{r,n}}{\sqrt{n}}\mcal{H}_\nu(n;h)=:\ts\sum_1+\cdots+\ts\sum_7
}
and treat $\ts\sum_\nu$'s in separate cases. \\

\subsection{Case $\ts\sum_1$} We claim that
$$
\ts\sum_1\ll T^3M.
$$
From \cite[(3.3.37) \& (3.3.38)]{Motohashi1997}
$$
(\widehat{h})'(\thf)=2\intt \f{\Gamma'(\thf+it)}{\Gamma(\thf+it)}t\pt h(t)\dd t \qand (\widehat{h})''(\thf)=4\intt \bigg\{\f{\Gamma'(\thf+it)}{\Gamma(\thf+it)}\bigg\}^{\!2}t\pt h(t)\dd t,
$$
we see
\aligs{
&(\widehat{h})'(\thf)=2i\pi^{\f{3}{2}}T^3M+O(TM^3)\\
\intertext{and}
&(\widehat{h})''(\thf)=8i\pi^{\f{3}{2}}T^3M\log T+O(TM^3\log T),
}
and thus
\aligs{
\ts\sum_1&=\f{2}{\pi^{3/2}}\big(T^3M\log T+O(T^3M)\big)\sum_{r\geq1}\f{1}{r}\sum_{n\geq1}\f{A_{r,n}d(n)}{n}\\
&\heq+\f{2}{\pi^{3/2}}\big(T^3M+O(TM^3)\big)\sum_{r\geq1}\f{1}{r}\sum_{n\geq1}\f{-A_{r,n}d(n)\log n}{n}.
}
Thus we need to establish
\alig{
\label{diag1}&\sum_{r\geq1}\f{1}{r}\sum_{n\geq1}\f{A_{r,n}d(n)}{n}\ll \f{1}{\log T}\\
\intertext{and}
\label{diag2}&\sum_{r\geq1}\f{1}{r}\sum_{n\geq1}\f{-A_{r,n}d(n)\log n}{n}\ll 1
}
In view of the identity
\aligs{
d(n_1n_2)=\sum_{s\mid(n_1,n_2)}\mu(s)d\Big(\f{n_1}{s}\Big)d\Big(\f{n_1}{s}\Big),
}
and that $a_n=0$ for $n\geq\xi^2$, we need for (\ref{diag1}) the bound
\alig{
\sum_{r\leq \xi^2}\f{1}{r}\sum_{s\leq \xi^2}\f{\mu(s)}{s^2}\sum_{n_1,n_2}\f{d(n_1)a_{rsn_1}\mu(rsn_1)}{n_1}\f{d(n_2)a_{rsn_2}\mu(rsn_2)}{n_2}\ll\f{1}{
\log T},
}
and for (\ref{diag2}) the bound
\alig{
\sum_{r\leq \xi^2}\f{1}{r}\sum_{s\leq \xi^2}\f{\mu(s)}{s^2}\sum_{n_1,n_2}\f{d(n_1)a_{rsn_1}\mu(rsn_1)\log n_1}{n_1}\f{d(n_2)a_{rsn_2}\mu(rsn_2)}{n_2}\ll 1.
}
But these last two bounds do hold, according to the same argument as in \cite[Section 2]{Luo2015}. \\

\subsection{Case $\ts\sum_2$} We write
\aligs{
\ts\sum_2&=\f{1}{2\pi^3}\sum_{r\leq\xi^2}\f{1}{r}\sum_{n\leq\xi^2}\f{A_{r,n}}{\sqrt{n}}\sum_{m=1}^{n-1}\f{d(m)d(m+n)}{\sqrt{m}}\Psi^+\Big(\f{m}{n};h\Big)\\
&\heq+\f{1}{2\pi^3}\sum_{r\leq\xi^2}\f{1}{r}\sum_{n\leq\xi^2}\f{A_{r,n}}{\sqrt{n}}\sum_{m\geq n}\f{d(m)d(m+n)}{\sqrt{m}}\Psi^+\Big(\f{m}{n};h\Big)\\
&=:\ts\sum_{2,1}+\ts\sum_{2,2}.
}
We have the bound which hold uniformly for $x\geq1$ (see \cite[p.\,123]{Motohashi1997}):
\alig{
\label{psi_p1}&\Psi^+(x;h)\ll x^{-1}T^{-B} \mbox{ for any fixed } 0<B<\sqrt{\log T}.
}
This bound implies that $\ts\sum_{2,2}$ is negligible.

By \cite[(3.4.20)]{Motohashi1997}, there is an absolute constant $c>0$ such that for $m<n$
\aligs{
\Psi^+\Big(\f{m}{n};h\Big)\ll T^3 M\exp\!\Big(\!-cM^2\f{m}{n}\Big)+\f{n}{m}T^{-\log T}.
}
For the range $T^\delta\log T\leq M<T/\log T$,  we have $M^2\f{m}{n}\geq(\log T)^2$ and thus
$$
\ts\sum_{2,1}\ll T^{4+2\delta-c\log T}+T^{4\delta-\qt\log T},
$$
which is negligible.

Next we consider the range $\log T<M<T^\delta\log T$ and follow Ivi\'{c} \cite{Ivic2001}. The nontrivial contribution to $\ts\sum_{2,1}$ comes from the sum
\aligs{
{\ts\sum^*_{2,1}}=\f{1}{2\pi^3}\sum_{r\leq\xi^2}\f{1}{r}\sum_{\f{M^2}{(\log T)^2}\leq n\leq\xi^2}\f{A_{r,n}}{\sqrt{n}}
\sum_{m\leq n \f{(\log T)^2}{M^2}}\f{d(m)d(m+n)}{\sqrt{m}}\Psi^+\Big(\f{m}{n};h\Big).
}
We abuse the notation $x=\f{m}{n}$, which is $o(1)$ as $T\rarrow\infty$. By \cite[(2.13)]{Ivic2001}, we have
\aligs{
\Psi^+(x;h)&=\f{4\pi\sqrt{x}}{\sqrt{x}+\sqrt{1+x}}\intt t\pt h(t)\tanh(\pi t)\\
&\heq\times\re\bigg\{\bigg(\!\f{\sqrt{x}+\sqrt{1+x}}{2}\bigg)^{\!\!-2it}\f{\Gamma^2(\hf+it)}{\Gamma(1+2it)}\\
&\heq\times F\bigg(\hf+it;\hf;1+it;\bigg(\!\f{\sqrt{x}-\sqrt{1+x}}{\sqrt{x}+\sqrt{1+x}}\bigg)^{\!\!2}\bigg)\bigg\}\dd t,
}
where $F(a,b;c;z)$ denotes the Gauss Hypergeometric function, initially defined for $|z|<1$ by
$$
F(a,b;c;z)=\sum_{k\geq0}\f{(a)_k(b)_k}{(c)_k}\f{z^k}{k!},
$$
where $(\alpha)_0=1$ and $(\alpha)_k=\alpha(\alpha+1)\cdots(\alpha+k-1)$ for $k\geq1$. Since
$$
\bigg(\!\f{\sqrt{x}-\sqrt{1+x}}{\sqrt{x}+\sqrt{1+x}}\bigg)^{\!\!2}=(\sqrt{x}+\sqrt{1+x})^{-4}<1-5\sqrt{x},
$$
we can use the absolute convergence of the hypergeometric series and write
\aligs{
\Psi^+(x;h)=\f{4\pi\sqrt{x}}{\sqrt{x}+\sqrt{1+x}}\sum_{k\geq0}\f{(\hf)_k}{k!}(\sqrt{x}+\sqrt{1+x})^{-4k}\re(J_k),
}
where
$$
J_k=\intt t\pt h(t)\tanh(\pi t)\f{(\hf+it)_k}{(1+it)_k}\bigg(\!\f{\sqrt{x}+\sqrt{1+x}}{2}\bigg)^{\!\!-2it}\f{\Gamma^2(\hf+it)}{\Gamma(1+2it)}\dd t.
$$

First we do some reduction. According to the concentration effect
of $e^{-\el(\!\f{t-T}{M}\!\er)^2}$ to $+T$ and $e^{-\el(\!\f{t+T}{M}\!\er)^2}$ to $-T$, we write
$$
J_k=J_k^++J_k^-
$$
and treat only $J_k^+$ since the two terms are very similar. For $J_k^+$, we further write
$$
J_k^+=\int_{|t-T|\leq M\log T}\cdots+\int_{|t-T|> M\log T}\cdots
$$
By Stirling's formula and that $|\Gamma(x+iy)|\geq \Gamma(x)\sqrt{\mathrm{sech}\pi y}$, we have
$$
\f{\Gamma^2(\hf+it)}{\Gamma(1+2it)}\ll (|t|+1)^{-\hf}\mbox{ for all real }t.
$$
Also $\el|\f{(\hf+it)_k}{(1+it)_k}\er|\leq1$. Thus it is easy to see that the second summand of $J_k^+$ contributes $\ll T^{\f{5}{2}-\log T}$, which can be neglected.
For the first summand of $J_k^+$, we bound trivially to obtain
$$
\int_{|t-T|\leq M\log T}\cdots\ll \int_{|t-T|\leq M\log T}t^{\f{5}{2}}\omega(t)\dd t\ll T^{\f{5}{2}}M,
$$
so that
$$
J_k^+\ll T^{\f{5}{2}}M.
$$
Hence, by the bound (see \cite[(2.17)]{Ivic2001})
\aligs{
\f{4\pi\sqrt{x}}{\sqrt{x}+\sqrt{1+x}}\sum_{k\geq0}\f{(\hf)_k}{k!}(\sqrt{x}+\sqrt{1+x})^{-4k}\ll x^{\qt},
}
we see that the contribution of $J_k^+$ to ${\ts\sum^*_{2,1}}$ is
\aligs{
&\ll T^{\f{5}{2}}M(\log T)^2\sum_{\f{M^2}{(\log T)^2}\leq n\leq\xi^2}\f{d(n)}{n^{\sss\f{3}{4}}}
\sum_{m\leq n \f{(\log T)^2}{M^2}}\f{d(m)d(m+n)}{m^{\sss\qt}}\\
&\ll T^{\f{5}{2}+2\delta}M,
}
which is $o(T^3M)$ upon letting $\delta<\f{1}{4}$.

We remark that finer analysis using the machinery in Ivi\'{c}'s work \cite{Ivic2001} leads to a larger range of $\delta$. This is only helpful if one could obtain an asymptotic formula
$$
\ts\sum_1=\mbox{some constant}\cdot T^3M+o(T^3M),
$$
which seems very difficult to achieve with our choice of the mollifiers $M_j$.\\

\subsection{Cases $\ts\sum_\nu$ ($\nu=3,\cdots,7$)}We shall see that the contribution of these $\ts\sum_\nu$'s can be neglected. \\

Case $\ts\sum_3$. The contribution of $\ts\sum_3$ is negligible due to the bound (\ref{psi_p1}).\\

Case $\ts\sum_4$. 
According to \cite[Section 4]{Ivic2001}, we have for $m<n$
\aligs{
\Psi^-\Big(\f{m}{n};h\Big)\ll TM^{\ep-1}\Big(\f{n}{m}\Big)^{\f{3}{2}-\ep}+T^{\ep-1}M\Big(\f{n}{m}\Big)^{\f{3}{2}}.
}
From this we see that
\aligs{
{\ts\sum_4}\ll T^{1+4\delta}\ll T^{2}
}
if we impose $\delta<\f{1}{4}$ as in the case of $\ts\sum_2$. \\

Case $\ts\sum_5$. The contribution of $\ts\sum_5$ is negligible because of the following bound (see \cite[p.\,123]{Motohashi1997})
\alig{
\label{psi_m1}&\Psi^-(x;h)\ll x^{-1}T^{-\log T},
}
which is uniform in $x\geq1$. \\

Case $\ts\sum_6$. The contribution of $\ts\sum_6$ is negligible since
\aligs{
{\ts\sum_6}=\f{12}{\pi^2}\el(e^{-\el(\!\f{T-i/2}{M}\!\er)^2}+e^{-\el(\!\f{T+i/2}{M}\!\er)^2}\er)\sum_{r\leq \xi^2}\f{1}{r}\sum_{n\leq\xi^2}\f{A_{r,n}}{\sqrt{n}}d_\hf(n).
}

Case $\ts\sum_7$. We simply discard $\ts\sum_7$ for its negativity, which is shown below. Define
$$
M_t=\sum_{n\geq1}\f{a_n\mu(n)d_{it}(n)}{\sqrt{n}}.
$$
Then
$$
|M_t|^2=\sum_{r\geq1}\f{1}{r}\sum_{n\geq1}\f{d_{it}(n)A_{r,n}}{\sqrt{n}},
$$
due to $d_{it}(m)d_{it}(n)=\sum_{d|(m,n)}d_{it}(mnd^{-2})$. Thus it follows from the definition of $\mcal{H}_7(n;h)$ that
\aligs{
\ts\sum_7&=-\f{1}{2\pi}\intt h(t)\f{|\zeta(\hf+it)|^4}{|\zeta(1+2it)|^2}\sum_{r\geq1}\f{1}{r}\sum_{n\geq1}\f{d_{it}(n)A_{r,n}}{\sqrt{n}}\dd t\\
&=-\f{1}{2\pi}\intt h(t)\f{|\zeta(\hf+it)|^4|M_t|^2}{|\zeta(1+2it)|^2}\dd t\\
&\leq0.
}

%
%
%
%
%
%
%
\section{A negative moment of $L(1,\sym^2u_j)$}\label{sec: neg_1st_moment_symmsq}
We closely follow Luo \cite{Luo2001} to prove Lemma \ref{negative_1st_moment_symmsq}, or the equivalent bound (see (\ref{hh_0}))
$$
\sum_j\f{h(t_j)}{L(1,\sym^2u_j)^2}\ll T^3M.
$$
For $\re(s)>1$, we write
\aligs{
\f{1}{L(s,\sym^2u_j)}=A_j(s)B_j(s),
}
where
$$
A_j(s)=\prod_p(1-\lambda_j(p^2)p^{-s})\qand B_j(s)=\prod_p\el(1+\f{\lambda_j(p^2)p^{-2s}-p^{-3s}}{1-\lambda_j(p^2)p^{-s}}\er).
$$
Since $B_j(s)$ is analytic and zero-free in $\re(s)\geq\f{9}{10}$, and $1/B_j(s)$ and $B_j(s)$ are uniformly bounded in this region, it suffices to show
\aligs{
\sum_j\f{A_j(1)}{L(1,\sym^2u_j)}h(t_j)\ll T^3M.
}
We have
$$
A_j(s)=\sum_{n\geq1}\f{\lambda_j(n^2)\mu(n)}{n^s}\quad(\re(s)>1).
$$
and by \cite[(36)]{Luo2001}
$$
A_j(1)=\sum_{n\geq1}\f{\lambda_j(n^2)\mu(n)e^{-n/T}}{n}-I_j
$$
where
$$
I_j=\f{1}{2\pi i}\int_{C_\ep}A_j(s+1)\Gamma(s)T^s\dd s
$$
and $C_\ep$ is the contour given by (\ref{contour}). Then we obtain
\aligs{
\sum_j\f{A_j(1)}{L(1,\sym^2u_j)}h(t_j)&=\sum_{n\geq1}\f{\mu(n)e^{-n/T}}{n}\sum_j\f{\lambda_j(n^2)}{L(1,\sym^2u_j)}h(t_j)-\sum_j\f{I_j}{L(1,\sym^2u_j)}h(t_j)\\
&=:S_1-S_2
}
For $S_1$, we apply Lemma \ref{KTF} to get
$$
S_1=D-C+O^+,
$$
where
\aligs{
D&=\hf\sum_{n\geq1}\f{\mu(n)e^{-n/T}}{n}\delta_{n,1}H,\\
C&=\hf\sum_{n\geq1}\f{\mu(n)e^{-n/T}}{n}\f{1}{\pi}\intt\f{d_{it}(n^2)}{|\zeta(1+2it)|^2}h(t)\dd t,\\
O^+&=\hf\sum_{n\geq1}\f{\mu(n)e^{-n/T}}{n}\sum_{c>0}\f{S(n^2,1;c)}{c}H^+\!\!\el(\f{4\pi n}{c}\er).
}
First it is easy to see that
$$
D=\hf e^{-1/T}H\ll \intt t\tanh(\pi t)h(t)\dd t\ll T^3M.
$$
Next we deduce that
\aligs{
C&\ll\sum_{n\geq1}\f{d(n^2)e^{-n/T}}{n}\intt\f{h(t)}{|\zeta(1+2it)|^2}\dd t\\
&\ll T^2M(\log T)^2\sum_{n\geq1}\f{d(n^2)e^{-n/T}}{n}\\
&\ll T^2M(\log T)^5
}
since
\aligs{
\sum_{n\geq1}\f{d(n^2)e^{-n/T}}{n}&=\f{1}{2\pi i}\int_{(2)}\sum_{n\geq1}\f{d(n^2)}{n^{1+s}}\Gamma(s)T^s\dd s\\
&=\f{1}{2\pi i}\int_{(2)}\f{\zeta(1+s)^3}{\zeta(2+2s)}\Gamma(s)T^s\dd s\\
&=\res{s=0}\f{\zeta(1+s)^3}{\zeta(2+2s)}\Gamma(s)T^s+\f{1}{2\pi i}\int_{C_\ep}T^s\Gamma(s)\f{\zeta(1+s)^3}{\zeta(2+2s)}\dd s\\
&\ll (\log T)^3,
}
where $C_\ep$ is the contour in (\ref{contour}). For
\aligs{
O^+\ll\sum_{n\geq1}\f{e^{-n/T}}{n}\sum_{c>0}\f{d(c)}{\sqrt{c}}\el|H^+\!\!\el(\f{4\pi n}{c}\er)\er|
}
we only need some control on $H^+(X)$ where we abuse the notation $X=\f{4\pi n}{c}$ for convenience. Shifting the integral of $H^+(X)$ to $\im(t)=-\sigma$ with $0<\sigma<\hf$, we get
\aligs{
H^+(X)&=\f{2i}{\pi}\intt\f{J_{2it}(X)}{\cosh(\pi t)}t\pt h(t)\dd t\\
&=\f{2i}{\pi}\int_{-\infty}^\infty \f{J_{2\sigma+2it}(X)}{\cosh(\pi (t-i\sigma))}(t-i\sigma)h(t-i\sigma)\dd t\\
&\ll T^{3-2\sigma}M X^{2\sigma},
}
due to (\ref{estJBessel}). Taking $\sigma=\qt+\f{\ep}{2}$, we have
\aligs{
O^+&\ll T^{3-\hf-\ep}M\sum_{n\geq1}\f{e^{-n/T}}{n^{\hf-\ep}}\sum_{c>0}\f{d(c)}{c^{1+\ep}}\\
&\ll T^{3-\hf-\ep}M\cdot\f{1}{2\pi i}\int_{(2)}\zeta(s+\thf-\ep)\Gamma(s)T^s\dd s\\
&=T^{3-\hf-\ep}M\bigg(\res{s=\hf+\ep}+\f{1}{2\pi i}\int_{\Gamma_\ep}\cdots\bigg)\\
&\ll T^3M,
}
where the contour $\Gamma_\ep=\el\{\hf+\ep+it\mid |t|\geq\ep\er\}\cup\el\{\hf+\ep+\ep e^{i\theta}\mid \f{\pi}{2}\leq\theta\leq\f{3\pi}{2}\er\}$. Summarizing the estimates of $D$, $C$ and $O^+$, we obtain that $S_1\ll T^3M$.

Now it remains to bound $S_2$. Let $J=\{j\mid |t_j-T|\leq M\log T\}$. We observe that
\aligs{
S_2\ll \sum_{j\in J}\f{I_j}{L(1,\sym^2u_j)}h(t_j),
}
since the contribution from $t_j<T-M\log T$ and $t_j>T+M\log T$ is negligible due to the bound
$$
I_j=A_j(1)+\f{1}{2\pi i}\int_{(\f{\ep}{2})}A_j(s+1)\Gamma(s)T^s\dd s\ll (t_jT)^\ep
$$
(see \cite[p.\,501]{Luo2001}).
For sufficiently small $\eta>0$, we partition $J$ into $J_1$ and $J_2$, according to whether $L(s,\sym^2u_j)$ is zero-free in $1-10\eta\leq\re(s)\leq1$ and $|\im(s)|\leq(\log T)^3$. By Luo's argument,
$$
I_j\ll T^{-\f{9}{20}\eta}\quad\mbox{for }j\in J_1.
$$
Hence by $L(1,\sym^2u_j)^{-1}\ll t_j^\ep$ (see \cite{HoffsteinLockhart1994}) and Weyl's law
\aligs{
&\sum_{j\in J_1}\f{I_j}{L(1,\sym^2u_j)}h(t_j)\\
&\ll T^{-\f{2}{5}\eta}\sum_{j\in J_1}h(t_j)\ll T^{2-\f{2}{5}\eta}\sum_{j\in J}1\\
&\ll T^{3-\f{2}{5}\eta}M\log T\\
&=o(T^3M).
}
On the other hand Luo's argument gives $|J_2|\ll T^{\f{1}{5}}$ and $\sum_{j\in J_2}|I_j|\ll T^{\qt}$. Thus
$$
\sum_{j\in J_2}\f{I_j}{L(1,\sym^2u_j)}h(t_j)\ll T^{2+\ep}\sum_{j\in J_2}|I_j|\ll T^{\f{5}{2}}.
$$
Hence we have shown that $S_2=o(T^3M)$ and completed the proof of Lemma \ref{negative_1st_moment_symmsq}.

%
%
%
%
%
%
%
\appendix
\section{A discussion on Barnes' formula}\label{sec: BarnesFormula}
In \cite{Barnes1900} Barnes developed the theory of the simple Gamma function $\Gamma_1(z\,|\,\omega)=\omega^{\f{z}{\omega}-1}\Gamma(z)$ with parameter $\omega$; $\Gamma_1(z\,|\,\omega)$ becomes $\Gamma(z)$ when $\omega=1$. Barnes's formula for $\Gamma(z)$ \cite[\S 41]{Barnes1900} states that for fixed $a\not\in\mbb{Z}\cap(-\infty,0]$ and for all large $z$ which are not in the vicinity of the negative real axis, we have
\alig{\label{Barnes}
\log\Gamma(z+a)=\Big(z+a-\hf\Big)\log z-z+{\hf}\log2\pi+\sum_{j=1}^n\f{(-1)^{j+1}B_{j+1}(a)}{j(j+1)z^n} + J_n(z;a),
}
where $n\geq0$ and $J_n(z;a)=O_a(|z|^{-n-1})$. Here $B_j(x)$ is the $j$-th Bernoulli polynomial, $\log z$ has the negative real axis as a cut and is real when $z$ is real and positive, and $J_n(z;a)$ is the contour integral (\ref{contourintegralJ}).
It is often useful to know the explicit dependence of the error term $J_n(z;a)$ on $a$, when dealing with ratios of Gamma functions. For this purpose, we prove the following
\prop{\label{BarnesError}Under the conditions for $a$ and $z$ in the above, we have
$$J_n(z;a)\ll P_{n+3}(|a|)|z|^{-n-1},$$
where the implied constant is absolute and $P_{n+3}(x)$ is a degree polynomial of degree $n+3$ whose coefficients are positive and may depend on $n$.
}
\begin{proof}
By the argument on \cite[p.\,121]{Barnes1900},
\eq{\label{contourintegralJ}
J_n(z;a)=-\int_{(\sigma)}\f{\pi}{s\sin\pi s}\,\zeta(s,a)z^s\dd s,
}
where $-n-1<\sigma<-n$ and $n\geq0$, and $\zeta(s,a)$ is the Hurwitz zeta function. Here we take $\sigma=-n-\hf$. By the argument in \cite[\S 40]{Barnes1900}, we can estimate $\zeta(s,a)$ with parameter $a$ as follows. For $-n-1<\re(s)=\sigma<-n$, we have
\aligs{
|\zeta(s,a)|\leq \sum_{m\leq[|a|]}\el|\f{1}{(m+a)^s}\er|+\f{C|a|^{n+2}}{(n+1)!}\sum_{m\geq[|a|]+1}\el|\f{1}{m^{s+n+2}}\er|+\sum_{m=0}^{n+1}\f{|a|^{m}|\zeta^{(m)}(s)|}{m!},
}
where $C>0$ is an absolute constant and $[\alpha]$ denotes the integer part of a real number $\alpha$. For $\sigma=-n-\hf$, we have
$$
\sum_{m\leq[|a|]}\el|\f{1}{(m+a)^s}\er|\leq \sum_{m\leq[|a|]}(m+|a|)^{n+1}
$$
and
$$
\f{C|a|^{n+2}}{(n+1)!}\sum_{m\geq[|a|]+1}\el|\f{1}{m^{s+n+2}}\er|\leq \f{C|a|^{n+2}}{(n+1)!}\sum_{m=1}^\infty\f{1}{m^{3/2}}.
$$
We also claim that the integrals
$$
\int_{(-n-\hf)}\el|\f{\pi}{s\sin\pi s}\,\zeta^{(m)}(s)\er||\dd s|
$$
are convergent for $0\leq m\leq n+1$ and will give the proof later. Collecting these estimates, we have
\aligs{
J_n(z;a)\ll Q_{n+2}(|a|)|z|^{-n-\hf},
}
where $Q_{n+2}(x)$ is a polynomial of degree $n+2$ with coefficients possibly dependent on $n$. With this bound to $J_{n+1}(z;a)$ and that
\aligs{
J_n(z;a)=\f{(-1)^{n+2}B_{n+2}(a)}{(n+1)(n+2)z^{n+1}}+J_{n+1}(z,a)\quad\mbox{(see \cite[p.\,120]{Barnes1900})},
}
the proposition follows.

Now we prove the claim. We have the functional equations
$$
\zeta(s)=2(2\pi)^{s-1}\sin\!\Big(\f{\pi s}{2}\Big)\Gamma(1-s)\zeta(1-s)
$$
and
$$
\zeta^{(m)}(s)=(-1)^m 2(2\pi)^{s-1}\sum_{j=0}^m\sum_{k=0}^m\el(a_{jmk}\sin\!\Big(\f{\pi s}{2}\Big)+b_{jmk}\cos\!\Big(\f{\pi s}{2}\Big)\er)\Gamma^{(j)}(1-s)\zeta^{(k)}(1-s),
$$
$m=1,\ldots,n+1$, as considered in \cite{Spira1965}, where $a_{jmk}$ and $b_{jmk}$ are constants independent of $s$. As shown in \cite{Spira1965}, for $s$ with $|s|\geq1$ and $|\arg (s)|<\pi$ one has
$$
\Gamma^{(j)}(s)=\Gamma(s)\el[(\log s)^j+\sum_{\ell=0}^{j-1}E_{\ell j}(s)(\log s)^\ell\er],
$$
where $E_{lj}(s)=O(|s|^{-1})$. Then, for $s=-n-\hf+it$,
\aligs{
\f{\pi}{s\sin\pi s}\,\zeta^{(m)}(s)&=\f{(-1)^m\pi(2\pi)^{s-1}}{s}
\sum_{j=0}^m\sum_{k=0}^m\el(\f{a_{jmk}}{\cos(\f{\pi s}{2})}+\f{b_{jmk}}{\sin(\f{\pi s}{2})}\er)\zeta^{(k)}(1-s)\\
&\heq\times\Gamma(1-s)\el[(\log(1-s))^j+\sum_{\ell=0}^{j-1}E_{\ell j}(1-s)(\log(1-s))^\ell\er]
}
is of exponential decay in $t$ as $t\rarrow\infty$. Thus we have shown the claimed convergence of the integrals
$$
\int_{(-n-\hf)}\el|\f{\pi}{s\sin\pi s}\,\zeta^{(m)}(s)\er||\dd s|,\quad m=1,\ldots, n+1.
$$
\end{proof}

\section{A uniform estimate for $K_{it}(x)$ ($t>0,\ x>1$)}\label{appendixB}
Here we give a simple consequence of the work of Booker--Str\"{o}mbergsson--Then \cite{BookerStrombergssonThen2013} on the $K$-Bessel function.
\prop{\label{Kit_estimate}For all $t>0$ and $x>1$,
\eq{\label{Kit_estimate_aux}
K_{it}(x)\ll e^{-\f{\pi}{2}t}t^{-\f{1}{3}}.
}
}
\pf{
First \cite[Proposition 2]{BookerStrombergssonThen2013} implies that for $x\geq t>0$
\aligs{
0<K_{it}(x)
&\leq e^{-\f{\pi}{2} t}e^{-\sqrt{x^2-t^2}+t\cos^{-1}(t/x)}
\min\el(\f{\sqrt{\f{\pi}{2}}}{(x^2-t^2)^{\f{1}{4}}},\f{\Gamma(\f{1}{3})}{2^{\f{2}{3}}3^{\f{1}{6}}}t^{-\f{1}{3}}\er)\\
&\ll e^{-\f{\pi}{2}t}t^{-\f{1}{3}}e^{-\sqrt{x^2-t^2}+t\cos^{-1}(t/x)}.
}
Then (\ref{Kit_estimate_aux}) holds for $x\geq t>0$, since
\eqs{
-\sqrt{x^2-t^2}+t\cos^{-1}(t/x)=t\el(\cos^{-1}(t/x)-\sqrt{(t/x)^{-2}-1}\er)\leq0
}
in view of
\eqs{
\f{\dd}{\dd u}\el(\cos^{-1}(u)-\sqrt{u^{-2}-1}\er)=\f{\sqrt{1-u^2}}{u^2}>0\quad\mbox{for }0<u<1.
}

On the other hand \cite[Proposition 2]{BookerStrombergssonThen2013} says that for $1\leq x<t$
\eqs{
|K_{it}(x)|<e^{-\f{\pi}{2}t}
\begin{cases}
\f{5}{(t^2-x^2)^{\f{1}{4}}}&\mbox{if }x\leq t-\hf t^{\f{1}{3}},\vspace*{1ex}\\
4t^{-\f{1}{3}}&\mbox{if }x\geq t-\hf t^{\f{1}{3}}.\end{cases}
}
For $x\leq t-\hf t^{\f{1}{3}}$, we have $t^2-x^2\geq t^{\f{4}{3}}-\f{1}{4}t^{\f{2}{3}}$. In addition, $t>1$ implies $t^{\f{4}{3}}-\f{1}{4}t^{\f{2}{3}}\geq\hf t^{\f{4}{3}}$. So (\ref{Kit_estimate}) holds when $1\leq x<t$.
}


\bibliographystyle{abbrv}


\end{document}